\documentclass[conference,letterpaper]{IEEEtran}

\usepackage[utf8]{inputenc} 
\usepackage[T1]{fontenc}
\usepackage{url}
\usepackage{ifthen}
\usepackage{cite}
\usepackage[cmex10]{amsmath} 
\usepackage{makecell}
\usepackage{caption} 
\usepackage{verbatim}
\usepackage{graphicx}
\usepackage{float}
\usepackage[colorinlistoftodos]{todonotes}
\usepackage{amssymb,amsfonts,bm}
\usepackage{amsthm}
\usepackage{tcolorbox}
\usepackage[framemethod=tikz]{mdframed}
\usepackage{mathtools,xparse}
\usepackage{tikz}
\usetikzlibrary{arrows,matrix,backgrounds}
\usetikzlibrary{shapes}
\usepackage{pgfplots}
\usetikzlibrary{plotmarks}
\pgfplotsset{compat=newest} 

\newtheorem{theorem}{Theorem}
\newtheorem{lemma}[theorem]{Lemma}

\def\cX{\mathcal{X}}
\def\cY{\mathcal{Y}}
\def\cP{\mathcal{P}}
\def\cQ{\mathcal{Q}}
\def\cD{\mathcal{D}}

\def\bp{\mathbf{p}}
\def\bpi{\mathbf{\pi}}

\def\pr{\text{Pr}}

\def\one{\mathbf{1}}
\def\GT{\widehat{M}_0^{\text{GT}}}

\begin{document}
\title{Missing Mass Estimation from Sticky Channels} 
\author{%
  \IEEEauthorblockN{Prafulla~Chandra, Andrew~Thangaraj}
  \IEEEauthorblockA{Department of Elec Engg\\
                    Indian Institute of Technology, Madras\\
                    Chennai, India 600036\\
                    \{ee16d402, andrew\}@ee.iitm.ac.in }
                    
                    \and 
        \IEEEauthorblockN{Nived~Rajaraman}
        \IEEEauthorblockA{Department of Elec Engg and Comp Sci \\
                    University of California, Berkeley\\
                    Berkeley, CA 94720, USA\\
                    nived.rajaraman@gmail.com}
        
}

\maketitle

\begin{abstract}
Distribution estimation under error-prone or non-ideal sampling modelled as "sticky" channels have been studied recently motivated by applications such as DNA computing. Missing mass, the sum of probabilities of missing letters, is an important quantity that plays a crucial role in distribution estimation, particularly in the large alphabet regime. In this work, we consider the problem of estimation of missing mass, which has been well-studied under independent and identically distributed (i.i.d) sampling, in the case when sampling is "sticky". Precisely, we consider the scenario where each sample from an unknown distribution gets repeated a geometrically-distributed number of times. We characterise the minimax rate of Mean Squared Error (MSE) of estimating missing mass from such sticky sampling channels. An upper bound on the minimax rate is obtained by bounding the risk of a modified Good-Turing estimator. We derive a matching lower bound on the minimax rate by extending the Le Cam method.
\end{abstract}

\section{Introduction}
Estimating distribution from samples is a fundamental question of interest in the fields of statistics, information theory and machine learning. This problem is particularly non-trivial in the large alphabet regime, i.e. when the number of samples is smaller or comparable to the alphabet size. The maximum likelihood (ML) estimator is not suitable for this regime since it assigns zero probability to the letters that are not seen in the samples, even if a considerable portion of the alphabet is not explored by the samples. To get a better estimate of the distribution in the large alphabet regime, we need to smooth the counts \cite{Chen1996} i.e. spread the counts of the letters seen in the samples on to the letters that were not seen in the samples. Smoothing, in general, requires an estimate of the  missing mass of the samples, which is the total probability of letters in the alphabet that were not seen in the samples. A popular estimator for missing mass is the Good-Turing (GT) estimator\cite{Good1953}, which is the ratio of the number of letters that occur only once in the samples to the sample size. The GT estimator and its variants are widely used for smoothing in natural language applications like speech recognition \cite{Chen1996}, \cite{Gale1995}, spelling correction \cite{Gale1991}, information retrieval \cite{Song99}, and ecology  \cite{Chao1992, Colwell12 }.

Missing mass has been studied theoretically when the samples generated are independent and identically distributed (i.i.d). Confidence intervals for missing mass were obtained using the GT estimator in \cite{McAllester2000} and were further improved in \cite{Chandra19}. The concentration of missing mass about its mean was studied in \cite{McAllester2003, berend2013, ben-hamou2017, Chandra19_ISIT}. The minimax risk of estimating missing mass of iid samples, from an unknown discrete alphabet, was first characterised in \cite{Rajaraman17} and was further improved in \cite{Acharya17}. The GT estimator and its several variations have been extensively studied and improved for the iid case \cite{ Acharya13, Orlitsky2015, Abs_disc_17, Ohannessian19, hao_doubly_comp19,  Skorski21,Painsky21}. In particular, the GT estimate for missing mass was used in \cite{Acharya13} for classification, pattern prediction, compression and in \cite{Orlitsky2015, hao_doubly_comp19} to estimate the underlying distribution on par with a natural genie estimator that knows the distribution. 

Recently there has been interest in estimating distributions or properties of distributions in the presence of sampling non-idealities. By sampling non-idealities, we refer to the smearing of data that could occur due to the practicalities in sampling a distribution (and making the samples non-iid) or the deletion or repetition of samples that generally occur due to lack of synchronization in data collection, storage and retrieval. For example, each sample can either be read multiple times or can not even be detected at all, during both storing and retrieval of data from DNA based storage devices \cite{Olgica_15}. In epidemiology, data about viral strains in a population is considered to be smeared since presence of a particular viral strain in an individual leads to the same strain being seen in many other individuals in the surroundings. These sampling non-idealities have been modelled in the literature as passing i.i.d samples from a distribution through a channel that repeats each sample $R$ times, where $R$ is a random variable on $\{0,1,\ldots\}$. Note that the sample gets deleted if $R = 0$. The channel is referred to as a sticky channel if $R$ is never 0 i.e. if no sample is deleted by the channel.

Estimating the support size of a distribution from the samples output from a Poisson repeat channel (where $R$ is Poisson) was studied in \cite{Olgica_21}. Distribution estimation using samples output from a geometric sticky channel, i.e. when $R$ follows a geometric distribution,  was studied in \cite{Olgica_09, Olgica_12}. Sticky channels have also been used to model sequencing errors in DNA sequencing \cite{nanopore_16}. The notion of geometric stickyness has been used in \cite{fried2021alphalazy} to make Markov trajectories ergodic and extend the results known for ergodic Markov chains. 

In this work, we consider the problem of estimating missing mass over a sticky channel. Specifically, we characterise the minimax rate of mean squared error (MSE) of estimating missing mass in samples output by a geometric sticky channel. Analysis of missing mass estimation is non-trivial even if the samples are i.i.d and becomes much tougher if the samples are output from a sticky channel that introduces memory in sampling.  

\subsection{Preliminaries}
We need the following definitions to set the context for our results. The geometric distribution, denoted Geo$(p)$, assigns probability $(1-p)^{k-1}p$ to $k\in\{1,2,\ldots\}$. Let $\bp \triangleq [p_1, \ldots, p_K]$ be a discrete distribution on the alphabet $\cX=[K]\triangleq\{1,2,\ldots,K\}$.  Let $C_S(\alpha)$ denote a geometric sticky channel that, for each input symbol, generates an independent $R\sim $ Geo$(1-\alpha)$ and repeats the input $R$ times.  
Let $X^n = (X_1,X_2,\ldots,X_n)$ be a sequence output by $C_S(\alpha)$ with iid samples from $\bp$ as the input.


For  $x \in \cX$, let 
$$N_x(X^n)\triangleq \sum_{i=1}^n I(X_i=x),$$ 
where $I(\cdot)$ is an indicator random variable, denote the number of occurrences of $x$ in $X^n$. 
The missing mass of $\bp$ in $X^n$ is the total probability of letters in $\cX$ that did not appear in $X^n$. Denoting the missing mass as $M_0(X^n,\bp)$, we see that
$$M_0(X^n,\bp) \triangleq \sum_{x \in \mathcal{X}} p_x I(N_x(X^n) = 0).$$

For $l\ge0$, let 
$$\phi_l(X^n) \triangleq \sum_{x\in\cX} I(N_x(X^n)=l)$$ 
denote the number of letters that have occurred $l$ times in $X^n$. The popular and standard Good-Turing estimator \cite{Good1953} for $M_0(X^n,\bp)$ is defined as
\begin{equation}
\GT(X^n) \triangleq \frac{\phi_1(X^n)}{n}. \label{eq: GT} \end{equation} 
We will drop the arguments $X^n, \bp$ whenever it is non-ambiguous.

We define the minimax risk, $ R_{n, \alpha}^{*},$  of estimating $M_{0}(X^n, \bp)$ from a sequence $X^n$ output by a geometric sticky channel $C_S(\alpha)$ with iid samples from $\bp$ as input, as 
\begin{align}
    R_{n,\alpha}^{*} \ & \triangleq \ \min_{ \widehat{M_0} } \ \max_{  \bp  } \  E [ (M_0(X^n, \bp) - \widehat{M_0}(X^n) )^2  ] \label{eq: minimax_gen}
\end{align}

\subsection{Markovian nature of $X^n$}
A Markov chain $Z^n$ over $\cX$ is the sequence $Z^n = (Z_1,Z_2,\ldots,Z_n)$ with the states $Z_i\in\cX$ and 
$$\pr(Z_{i}{=}z_i|Z_{i-1}{=}z_{i-1},\ldots,Z_1{=}z_1)=\pr(Z_2{=}z_i|Z_1{=}z_{i-1})$$ 
for $i=2,\ldots,n$ and all $z_i\in\cX$. The transition probability matrix (t.p.m) of the Markov chain, denoted $P$, is the $K\times K$ matrix with $(i,j)$-th element $P_{ij}$ (or $P(j|i))$ $\triangleq \pr(Z_2=j | Z_1 = i)$. A distribution $\bpi=[\pi_1,\ldots,\pi_K]$ on $\cX=[K]$ is said to be a stationary or invariant distribution of the Markov chain if $\bpi P = \bpi$ \cite{markov_book}. 

The geometric stickyness of $C_{S}(\alpha)$ leads to the sequence $X^n$ being generated as follows:
the channel $C_{S}(\alpha)$ draws a sample from $\bp$ and outputs it as $X_1.$ 
For $i \geq 2,$ the channel sets $X_{i}$  
\begin{enumerate}
    \item to $X_{i-1},$ with probability $\alpha,$  
    \item to an independent sample from $\bp,$ with probability $1-\alpha.$ 
\end{enumerate}
Therefore the sequence $X^n$ is a Markov chain with 
\begin{align}
& \pr(X_i = z|X_{i-1} = y) \ = \ 
\alpha I(z = y) + (1-\alpha) p_y,  \nonumber \\
& \text{ for } y,z \in \cX, \ i = 2,\ldots,n,  \label{eq: Xn_MC}
\end{align}
i.e. the sequence $X^n$ is a stationary Markov chain with t.p.m $P = \alpha I_{K \times K} + (1-\alpha) \ \one_{K \times 1} \ \bp$, where $\one_{K \times 1}$ is the all ones column vector of length $K$, and with $\bp$ as the stationary distribution. 
\subsection{Estimators for $M_0$ of sticky samples}
We first consider the bias of a scaled version of the Good-Turing estimator \eqref{eq: GT} in estimating $M_0$ from sticky samples. For $m = 1$ to $n$, let $X_{\sim m} \triangleq (X_1,\ldots, X_{m-1}, X_{m+1},\ldots X_n)$ denote the samples from $X_1$ to $X_n$ except $X_m$. Let $X_{i}^j \triangleq (X_i, \ldots, X_j)$ denote the samples from $X_i$ to $X_j$.
\begin{align*}
    & E\left[\kappa\GT(X^n) - M_0(X^n,\bp)\right] \\ 
    & = \ \frac{\kappa}{n} E[\phi_1] - \sum_{x \in \cX} p_x \ \pr(N_x = 0) \\
    & = \ \sum_{x \in \cX} \left(\frac{\kappa}{n} \pr(N_x = 1) - p_x \ \pr(N_x = 0)\right) \\
    & = \  \frac{1}{n} \sum_{m=1}^n \left(\sum_{x \in \cX} \kappa\pr(X_m = x, N_x(X_{\sim m}) = 0)  -  p_x \ \pr(N_x = 0)\right).
\end{align*}
For a sequence $X^n$ output by a geometric sticky channel $C_{S}(\alpha)$ with iid samples from $\bp$ as input, we have
\begin{enumerate}
    \item $\pr(N_x(X^n) = 0) \overset{(a)}{=} (1-p_x) \  (1-(1-\alpha)p_x)^{n-1},$ \\
    \item for $2\leq m \leq n-1,$ \\
    $ \pr(X_m = x, N_x(X_{\sim m}) = 0) \\
    \overset{(b)}{=} \ p_x \ (1-\alpha)^2 \ (1-p_x)^2 \ (1-(1-\alpha)p_x)^{n-3},$
\end{enumerate} 
where we get $(a),(b)$ by using \eqref{eq: Xn_MC}. So for m = $2$ to $n-1,$
 \begin{align}
     & \sum_{x \in \cX} \kappa\pr(X_m = x, N_x(X_{\sim m}) = 0) - p_x \ \pr(N_x = 0) \nonumber \\ 
   & = \ \sum_{x \in \cX} p_x(1-p_x)(1-(1-\alpha)p_x)^{n-3}\nonumber\\
   &\qquad\qquad [\kappa (1-\alpha)^2 (1-p_x) -  (1-(1-\alpha)p_x)^2] \nonumber
\end{align}
Setting $\kappa=(1-\alpha)^{-2}$ results in cancellation and the term above decays as $1/n$ and a vanishing bias. 

Motivated by the above bias calculation, we define the modified Good-Turing estimator for the missing mass of a sequence $X^n$ output by a geometric sticky channel $C_S(\alpha)$ as 
\begin{equation}
\GT(X^n,\widehat{\alpha}) \triangleq \frac{\phi'_1(X^n)}{(1-\widehat{\alpha})^2(n-2)}. \label{eq: GT_mod} \end{equation} 
where $\phi'_1(X^n) \triangleq \sum_{x\in\cX} I(N_x(X_2^{n-1})=1, X_1 \neq x, X_n \neq x)$ is the number of singletons of $X^n$ that have occurred in $X_2,\ldots,X_{n-1}$. We set $\widehat{\alpha}$ to  $\alpha$, if $\alpha$ is known, or to an estimate of $\alpha$, if $\alpha$ is unknown. 

The division by $(1-\alpha)^2$ is natural for two reasons: let $Y^t = (Y_1,\ldots, Y_t)$ ($t$ unknown), be the samples (iid $\sim \bp$) input to $C_S(\alpha)$ that resulted in $X^n.$ Since each $Y_i$ is repeated for $R \sim$ Geo($1-\alpha)$ times by $C_S(\alpha),$ on an average,
\begin{enumerate}
    \item  $(1-\alpha)$ of the singletons in $Y^t$ remain singletons in $X^n,$ i.e. $\phi_1(X^n)/(1-\alpha)$ is a good estimate of the number of singletons in $Y^t,$ 
    \item  number of repetitions per input is $1/(1-\alpha),$ i.e. $(1-\alpha) n$ is a good estimate of $t.$
\end{enumerate}
If $Y^t = (Y_1,\ldots, Y_t)$ were known to us, the Good-Turing estimate for $M_0(Y^t,\bp)$ would be $ \phi_1(Y^t)/t.$ Since $M_0(X^n,\bp) = M_0(Y^t,\bp),$
 we get the modified Good-Turing estimator in \eqref{eq: GT_mod} by using $\phi_1(X^n)/(1-\alpha)$ and $(1-\alpha) n$ in the place of $\phi_1(Y^t)$ and $t$, respectively.\footnote{We omit $X_1, X_n$ and consider only $X_2,\ldots,X_{n-1}$ to make the analysis simpler.} 

\section{Results}
Our main result is the characterisation of the minimax rate $R^*_{n,\alpha}$ through upper and lower bounds that match in the order of $n$.  The first result is an upper bound on $R_{n,\alpha}^{*}$ obtained by characterising the worst case MSE of the modified Good-Turing estimator $\GT(X^n,\alpha),$ assuming $\alpha$ is known. 
\begin{theorem}
\begin{align}
    R_{n,\alpha}^{*} \ & \leq \ \frac{1}{(n{-}2)(1{-}\alpha) } \left( 3 + \frac{1}{1{-}\alpha} \left(1 + \frac{1}{n}\right) \right) + O(1/n)  \label{eq: minimax_ub} 
\end{align}
\label{thm: minimax_ub}
\end{theorem}
\begin{proof}
Section \ref{sec: minimax_ub_Prf}.
\end{proof}
For $\alpha = 0$ i.e if $X^n \sim \bp$ iid, the bound in \eqref{eq: minimax_ub} reduces to the $O(1/n)$ bound on the worst case MSE of the Good-Turing estimator $\GT(X^n).$
\begin{theorem}
\begin{align}
    R_{n,\alpha}^{*} \ & \geq \frac{1/32}{1 + (n{-}1) (1{-}\alpha)} - \left( 1 - 0.5(1{-}\alpha) \right)^{n-1}, \label{eq: minimax_lb}
\end{align}
\label{thm: minimax_lb}
\end{theorem}
\begin{proof}
Section \ref{sec: minimax_lb_Prf}
\end{proof}
We obtain the lower bound in \eqref{eq: minimax_lb} by modifying the classical Le Cam method. 

If $\alpha$ is unknown, we consider the estimate $\widehat{\alpha}=1 - (\tau(X^n) + 1)/n$, where $\tau(X^n)$ is the number of state changes in $X^n$. Through simulations, we observe that $\widehat{\alpha}$ is a good estimate if $\alpha \gg \max_{x \in \cX} p_x$, which  can be expected to hold for many $C_S(\alpha), \bp$ pairs in the large alphabet regime.  

Fig. \ref{fig:MSE_plaw_0.1}, \ref{fig:MSE_plaw_0.5_0.1}  plot the MSE of $\GT(X^n, \widehat{\alpha})$ against $n$ with the i.i.d samples from 
\begin{enumerate}
    \item the power law distribution having $p_i \propto 1/i^{0.1}$ on  $\cX = \{1,\ldots,K\}, K = 1.2n,$ 
    \item the nearly power law distribution with $p_1 = 0.1,$ $p_i \propto 1/i^{0.5}$ on  $\cX = \{2,\ldots,K\}, K = 1.2n,$ 
\end{enumerate}
respectively as the input to $C_{S}(\alpha),$ for $\alpha = 0.5,0.75,0.95$.
\begin{figure}[h]
\begin{center}
%
%
\definecolor{mycolor1}{rgb}{1.00000,0.00000,1.00000}%
\begin{tikzpicture}

\begin{axis}[%
width=2.5in,
at={(0.758in,0.481in)},
scale only axis,
xmode=log,
xmin=100,
xmax=10000,
xtick={100,400,1600,6400},
xticklabels={{100},{400},{1600},{6400}},
xminorticks=true,
xlabel style={font=\color{white!15!black}},
xlabel={sample size n},
ymode=log,
ymin=0.0001,
ymax=10,
yminorticks=true,
ylabel style={font=\color{white!15!black}},
ylabel={Mean Squared Error of modified GT},
axis background/.style={fill=white},
title style={font=\bfseries},
title={Power law $p_i \propto 1/i^{0.1}$ on $\{1,\ldots, K\}$},
xmajorgrids,
xminorgrids,
ymajorgrids,
yminorgrids
]
\addplot [color=black, mark=asterisk, mark options={solid, black}, forget plot]
  table[row sep=crcr]{%
100	0.0370428364234995\\
200	0.0177274149940151\\
400	0.00903722770663386\\
800	0.00448136368643165\\
1600	0.0022704899784561\\
3200	0.00110766721401108\\
6400	0.000561771303009398\\
};
\addplot [color=blue, dashed, mark=x, mark options={solid, blue}, forget plot]
  table[row sep=crcr]{%
100	0.0194\\
200	0.0095\\
400	0.0047\\
800	0.0023\\
1600	0.0012\\
3200	0.0006\\
6400	0.0003\\
};
\addplot [color=mycolor1, mark=diamond, mark options={solid, mycolor1}, forget plot]
  table[row sep=crcr]{%
100	0.17751795950763\\
200	0.0864241042153257\\
400	0.0426401176878221\\
800	0.021568494733063\\
1600	0.0105438660062992\\
3200	0.00530360457983361\\
6400	0.00264789386656172\\
};
\addplot [color=red, dashdotted, mark=star, mark options={solid, red}, forget plot]
  table[row sep=crcr]{%
100	0.0962\\
200	0.0459\\
400	0.0231\\
800	0.0114\\
1600	0.0058\\
3200	0.0029\\
6400	0.0014\\
};
\addplot [color=black, mark=triangle, mark options={solid, rotate=270, black}, forget plot]
  table[row sep=crcr]{%
100	5.8776\\
200	2.4314\\
400	1.1344\\
800	0.5482\\
1600	0.2686\\
3200	0.1307\\
6400	0.0654\\
};
\addplot [color=blue, dashed, mark=square, mark options={solid, blue}, forget plot]
  table[row sep=crcr]{%
100	2.9243\\
200	1.7702\\
400	0.8794\\
800	0.4419\\
1600	0.2187\\
3200	0.1088\\
6400	0.0546\\
};
\end{axis}

\begin{axis}[%
width=2.5in,
at={(0in,0in)},
scale only axis,
xmin=0,
xmax=1,
ymin=0,
ymax=1,
axis line style={draw=none},
ticks=none,
axis x line*=bottom,
axis y line*=left
]
\node[below right, align=left, draw=white]
at (rel axis cs:0.65,0.45) {$\alpha\text{ = 0.5}$};
\node[below right, align=left, draw=white]
at (rel axis cs:0.77,0.79) {$\alpha\text{ = 0.75}$};
\node[below right, align=left, draw=white]
at (rel axis cs:0.35,1.02) {$\alpha\text{ = 0.95}$};
\end{axis}

\end{tikzpicture}%
\end{center}
\captionof{figure}{MSE v/s $n$ for power law distribution.}
\label{fig:MSE_plaw_0.1}
\end{figure}
\begin{figure}[h]
\begin{center}
%
%
\definecolor{mycolor1}{rgb}{1.00000,0.00000,1.00000}%
\begin{tikzpicture}

\begin{axis}[%
width=2.5in,
at={(0.758in,0.481in)},
scale only axis,
xmode=log,
xmin=100,
xmax=10000,
xtick={100,400,1600,6400},
xticklabels={{100},{400},{1600},{6400}},
xminorticks=true,
xlabel style={font=\color{white!15!black}},
xlabel={sample size n},
ymode=log,
ymin=0.0001,
ymax=10,
yminorticks=true,
ylabel style={font=\color{white!15!black}},
ylabel={Mean Squared Error of modified GT},
axis background/.style={fill=white},
title style={font=\bfseries},
title={Nearly Power law $p_1 = 0.1, p_i \propto 1/i^{0.5}$ on $\{2,\ldots, K\}$},
xmajorgrids,
xminorgrids,
ymajorgrids,
yminorgrids
]
\addplot [color=black, mark=asterisk, mark options={solid, black}, forget plot]
  table[row sep=crcr]{%
100	0.0307765375183394\\
200	0.0145594439143963\\
400	0.00723528220063437\\
800	0.00357075536638047\\
1600	0.00185268876602271\\
3200	0.000908310509218837\\
6400	0.00044863\\
};
\addplot [color=blue, dashed, mark=x, mark options={solid, blue}, forget plot]
  table[row sep=crcr]{%
100	0.0196180583371804\\
200	0.00948466388065091\\
400	0.00467048888880328\\
800	0.00241334517395176\\
1600	0.00125814464071685\\
3200	0.000692291073054999\\
6400	0.00039813\\
};
\addplot [color=mycolor1, mark=diamond, mark options={solid, mycolor1}, forget plot]
  table[row sep=crcr]{%
100	0.149539384538093\\
200	0.0737784071991\\
400	0.0356311043704992\\
800	0.0171400769908473\\
1600	0.00847600146152074\\
3200	0.00427032203534144\\
6400	0.00210870059561642\\
};
\addplot [color=red, dashdotted, mark=star, mark options={solid, red}, forget plot]
  table[row sep=crcr]{%
100	0.0924209500063942\\
200	0.0454165907915669\\
400	0.0223763439852623\\
800	0.0108535101751669\\
1600	0.00549719766782192\\
3200	0.00282608182209013\\
6400	0.00152284316735179\\
};
\addplot [color=black, mark=triangle, mark options={solid, rotate=270, black}, forget plot]
  table[row sep=crcr]{%
100	5.30793675139182\\
200	2.2146582920055\\
400	1.04292251057106\\
800	0.48173247410446\\
1600	0.235775494037058\\
3200	0.112053167746143\\
6400	0.0571622544291893\\
};
\addplot [color=blue, dashed, mark=square, mark options={solid, blue}, forget plot]
  table[row sep=crcr]{%
100	3.05860254655735\\
200	1.72927368169386\\
400	0.855366368856024\\
800	0.415027830582552\\
1600	0.204543394914871\\
3200	0.0995732355895942\\
6400	0.0506098445119165\\
};
\end{axis}

\begin{axis}[%
width=2.5in,
at={(0in,0in)},
scale only axis,
xmin=0,
xmax=1,
ymin=0,
ymax=1,
axis line style={draw=none},
ticks=none,
axis x line*=bottom,
axis y line*=left
]
\node[below right, align=left, draw=white]
at (rel axis cs:0.40,0.55) {$\alpha\text{ = 0.5}$};
\node[below right, align=left, draw=white]
at (rel axis cs:0.78,0.75) {$\alpha\text{ = 0.75}$};
\node[below right, align=left, draw=white]
at (rel axis cs:0.33,1.02) {$\alpha\text{ = 0.95}$};
\end{axis}

\end{tikzpicture}%
\end{center}
\captionof{figure}{MSE v/s $n$ for a nearly power law distribution.}
\label{fig:MSE_plaw_0.5_0.1}
\end{figure}
We performed 16000 trials for $n \in  \{100,200,400,800,1600,3200,6400\}.$ 
The solid lines in Fig. \ref{fig:MSE_plaw_0.1}, \ref{fig:MSE_plaw_0.5_0.1} is the MSE of $\GT(X^n,\alpha)$ i.e. assuming $\alpha$ to be known. The dashed lines represent the MSE of $\GT(X^n,\widehat{\alpha})$. The success of $ 1-(\tau(X^n) + 1)/n$ in estimating $\alpha$ is intuitive because $\alpha \gg \max_{x \in \cX} p_x$ for the power law distributions considered above leading to $\max_{x \in \cX} P_{ii} = \alpha + (1-\alpha) p_i \approx \alpha$.

The rest of the paper is devoted to the proofs of Theorems \ref{thm: minimax_lb} and \ref{thm: minimax_ub}.

\section{Proof of Theorem \ref{thm: minimax_ub} }
\label{sec: minimax_ub_Prf}
We get the upper bound in \eqref{eq: minimax_ub} by bounding the worst case MSE of the modified Good-Turing estimator in \eqref{eq: GT_mod}. The MSE of the modified Good-Turing estimator is first expressed as follows. 
\begin{lemma}
For a sequence $X^n$ output from $C_S(\alpha)$ with i.i.d samples from $\bp$ as input, the MSE of the modified Good-Turing estimator $\GT(X^n, \alpha)$ is given by 
\begin{align}
& E[(\GT(X^n,{\alpha}) - M_0(X^n,\bp))^2] \ = \  \sum_{x,y \in \cX : x \neq y} T_{x,y}(n,\alpha) \nonumber \\ 
& \quad + \sum_{x \in \cX} p_x^2 \ Q_x(0) + \frac{1}{(1-\alpha)^4 (n-2)^2} \sum_{x \in \cX} Q_x(1),  \label{eq: MSE_expr}
\end{align}
where $Q_x(a)\triangleq \pr(N^{\sim 1,2}_x = a, E^{\sim x}_{1,n})$ with the notation $N^{\sim 1,2}_x=N_x(X_2^{n-1})$ and $E^{\sim x}_{1,n} \triangleq (X_i \neq x, i = 1,n)$, and
\begin{align}
Q_{x,y}(a,b) &\triangleq \pr(N^{\sim 1,2}_x = a, N^{\sim 1,2}_y = b, E^{\sim x}_{1,n},E^{\sim y}_{1,n}), \nonumber \\
T_{x,y}(n,\alpha) &\triangleq p_x p_y  Q_{x,y}(0,0) + \frac{1}{(1-\alpha)^4 (n-2)^2}Q_{x,y}(1,1)\nonumber \\
& \quad - \frac{1}{(1-\alpha)^2 (n-2)}(p_x Q_{x,y}(0,1) + p_y Q_{x,y}(1,0)),\nonumber
 \end{align}
\label{lem: MSE_expr}
\end{lemma}
\begin{proof}
Section \ref{subsec: MSE_expr_Prf} in Appendix. 
\end{proof}
The following lemma bounds each of the terms in \eqref{eq: MSE_expr}.
\begin{lemma}
\begin{align}
&\sum_{x \in \cX} p_x^2 \ Q_x(0)  \  \leq \ (1-\alpha)^{-1}/(n+1), \label{eq: t2_bd} \\
&\frac{1}{(1-\alpha)^4 (n-2)^2}  \sum_{x \in \cX} Q_x(1) \  \leq \ (1-\alpha)^{-2}/(n-2), \label{eq: t3_bd} \\
& \sum_{x,y \in \cX: x \neq y} T_{x,y}(n,\alpha) \  \leq \ 2(1-\alpha)^{-1}/(n-2) + O(1/n) \nonumber \\ 
& \qquad + \left( \alpha^2 \left( \frac{\alpha^2}{(1-\alpha)^2} + 2 \right) +  (1-\alpha)^2 \right) O(1/n^2). \label{eq: t1_bd}  
\end{align}
\label{lem: MSE_term_bds}
\end{lemma}
\begin{proof}
Section \ref{sec: MSE_term_bds_Prf} in Appendix. 
\end{proof}

Using the bounds in \eqref{eq: t2_bd}, \eqref{eq: t3_bd}, and \eqref{eq: t1_bd} in \eqref{eq: MSE_expr} completes the proof of \eqref{eq: minimax_ub}.

\section{ Proof of Theorem \ref{thm: minimax_lb}  } 
\label{sec: minimax_lb_Prf}
The standard Le Cam method \cite{yu1997lecam} is for estimating constant parameters of a distribution whereas $M_0$ is a function of both the distribution and the samples. To prove Theorem \ref{thm: minimax_lb}, we first extend the Le Cam method to the context of estimating functions of both the distribution and its samples.

\subsection{Le Cam lower bound for estimating random variables}
\label{subsec:Lecam_lb}
Let $ \cQ $ be a family of distributions over an alphabet $\cY$ and $Y$ be a random variable distributed according to $Q \in \cQ$. Let $\theta(Y,Q)$, taking values in a pseudometric space $\cD$ with a pseudometric $d$, be a function of both $Y$ and the distribution $Q$. We assume that the set $\cD$ is bounded i.e. the distance $d(u,v)$ between any two points $u , v \in \cD$ is at most $\Delta$.  
Let $d(\cD_1, \cD_2) \triangleq \min_{u \in \cD_1, v \in \cD_2} d(u,v)$ be the distance between the subsets $\cD_1, \cD_2$ of $\cD$. Let $\widehat{\theta}(Y)$ be an estimator for $\theta(Y,Q)$ and $co(\cQ)$ denote the convex hull of $\cQ$. 

The next lemma provides a lower bound on the worst case risk (over $\cQ$) of any estimator $\widehat{\theta}(Y)$ for $\theta(Y,Q)$. 
\begin{lemma}
Let $ \cD_1, \cD_2$ be two subsets of $\cD$, and  $\cQ_1, \cQ_2$ be two subsets of $\cQ$ such that for any $Q_i \in \cQ_i,$ $\theta(Y,Q_i) \in \cD_i$ with probability at least $1 - \epsilon_i$, $i = 1,2$. Let  $\delta \triangleq d(\cD_1, \cD_2)/2,$ and $||Q_1 \wedge Q_2 || \triangleq 1 - ||Q_1 - Q_2||_{TV}$ denote the affinity of the two distributions $Q_1$ and $Q_2$. Then 
\begin{align}
   & \sup_{Q \in \cQ} E [ d( \widehat{\theta}(Y) , \theta(Y,Q) ) ] \nonumber \\ 
   & \ \geq \ \delta \ \left( \sup_{Q_i \in co(\cQ_i)} ||Q_1 \wedge Q_2 || \right) - (\epsilon_1 + \epsilon_2)\Delta.  \label{eq:Lecam_lb} 
\end{align}
\label{lem:Lecam_lb}
\end{lemma}
\begin{proof}
Section \ref{subsec:Lecam_lb_Prf} in Appendix.
\end{proof}
To use Lemma \ref{lem:Lecam_lb} in the context of sticky channels, we let $\cQ=\{P_{\bp}(X^n):X^n\text{ output of }C_S(\alpha) \text{ with input i.i.d }\bp\}$ be the space of all distributions from sticky channels. The output $X^n$ is considered as the random variable $Y$ and $\Theta=M_0(X^n,\bp)$. 

Consider the class of distributions $\bp(\gamma,L) = \begin{bmatrix} 0.5+\gamma & ((0.5-\gamma)/{L})\one_{1 \times L}  \end{bmatrix}^T$ on the alphabet $\{1,\ldots,L+1\}$, where $\one_{1 \times L}$ denotes the length-$L$ all-1s vector. In $X^n$ obtained using input distribution $\bp(\gamma,L)$, the observed mass is maximized when 1 occurs along with $n-1$ other distinct letters. So, 
$$M_{\text{obs}}\triangleq1-M_0(X^n,\bp(\gamma,L))\le 0.5+\gamma+n(0.5-\gamma)/L.$$
If 1 occurs at least once, the observed mass is $\ge 0.5+\gamma$. So,
$$M_{\text{obs}}\ge 0.5+\gamma \text{ w. p. }\ge 1-(0.5-\gamma)(1-(1-\alpha)(0.5+\gamma))^{n-1}.$$
Combining, we get that 
$$M_0(X^n,\bp(\gamma,L))\in(0.5-\gamma)[1-n/L,1]$$ 
with high probability.

Using the above observation, we proceed to let the subsets of distributions $\cQ_1$ and $\cQ_2$ needed in Lemma \ref{lem:Lecam_lb} as $\cQ_1=P_{\bp(0,L)}(X^n)$ and $\cQ_2=P_{\bp(\beta,L)}(X^n)$ for some $\beta$, $n/2L<\beta<0.5$, and $L$ a large integer to be chosen later. The subset $\cD_1=0.5[1-n/L,1]$ contains the missing mass $M_0(X^n,\bp(0,L))$ w.p. $1-\epsilon_1$ for $\epsilon_1=0.5(1-(1-\alpha)0.5)^{n-1}$. Similarly, the subset $\cD_2=(0.5-\beta)[1-n/L,1]$ contains $M_0(X^n,\bp(\beta,L))$ w.p. $1-\epsilon_2$ for $\epsilon_2=(0.5-\beta)(1-(1-\alpha)(0.5+\beta))^{n-1}$. The sets $\cD_1$ and $\cD_2$ are separated by a distance $\beta-n/2L$.

Plugging into Lemma \ref{lem:Lecam_lb} with $d(u,v)=(u-v)^2$, 
\begin{align}
     & \sup_{\bp \in \cP} E [ ( \widehat{M_0}(X^n) - M_0(X^n,\bp) )^2 ] \nonumber \\
     & \ \geq 0.5  \left( {\beta} - \frac{n}{2L} \right)^2  \lVert P_{\bp(0,L)} \wedge P_{\bp(\beta,L)}\rVert  - (\epsilon_1+\epsilon_2). \label{eq:M0_Lecam_int}
\end{align}
The next lemma provides a simplified bound on the total variation distance between $P_{\bp(0,L)}(X^n)$ and $P_{\bp(\beta,L)}(X^n).$
\begin{lemma}
\begin{equation}
    \lVert P_{\bp(0,L)} - P_{\bp(\beta,L)}\rVert_{TV} \leq \sqrt{2} {\beta} \left(  1 + 0.5  (n-1) (1-\alpha) \right)^{0.5}. \label{eq:TV_bd} 
\end{equation}
\label{lem:TV_bd}
\end{lemma}
\begin{proof}
Section \ref{subsec:TV_bd_Prf} in Appendix.
\end{proof}
Using in $\lVert P_{\bp(0,L)} \wedge P_{\bp(\beta,L)}\rVert$, we get 
\begin{align}
    & \sup_{\bp \in \cP} E [ ( \widehat{M_0}(X^n) - M_0(X^n,\bp) )^2 ] \nonumber \\ 
    & \geq 0.5 \left( \beta - \frac{n}{2L} \right)^2 \left( 1 - \sqrt{2} \beta \left( { 1 + 0.5  (n-1) (1-\alpha)  } \right)^{0.5}  \right) \nonumber \\ 
    & \qquad -\left( 1 - 0.5(1-\alpha) \right)^{n-1}. \label{eq:M0_Lecam_int2}
\end{align}
Setting $\beta = \frac{1}{2 \sqrt{2}} \ \left( {  1 + 0.5  (n-1)  (1-\alpha)   } \right)^{-0.5},$ and  $L =e^n$ in \eqref{eq:M0_Lecam_int2} completes the proof of Theorem \ref{thm: minimax_lb}.

\section{Conclusion and Future Directions}
Based on our analytical study, we conclude that the problem of estimation of missing mass is tractable in the presence of sampling non-idealities such as those that can be modelled as outputs of sticky channels. Since missing mass (with ideal sampling) is generally considered to be a non-trivial problem with a surprisingly elegant estimator, it is all the more surprising that it remains robust to non-idealities such as sticky channels that introduce memory. A simple modification to the Good-Turing estimator is seen to be minimax-optimal when the stickiness probability is known. An extension of the Le Cam method provides the matching lower bound.

A natural future direction would be to characterise the minimax risk of estimating combined probability mass, i.e. total probability of letters that have each occurred $k ( \geq 1)$ times, in the samples output by a geometric sticky channel. Another extension would be to characterise the minimax risk of estimating missing mass of samples output by other channels such as a duplication channel i.e a sticky channel in which the number of repetitions is a Bernoulli random variable. 

\appendix

\subsection{Proof of Lemma \ref{lem: MSE_expr}}
\label{subsec: MSE_expr_Prf}
Substituting  $M_{0}(X^n, \bp) $ $ =  \sum_{x \in \cX} p_x I(N_x(X^n) = 0)$ and $\GT(X^n,\alpha) = \frac{1}{(1-\alpha)^2 (n-2) }  \sum_{x \in \cX} I(N_x(X_2^{n-1}) = 1, E^{\sim x}_{1,n})$ into $E[(M_0(X^n,\bp) - \GT(X^n,\alpha))^2] $ and taking expectation of each term in the square of the summation, we have
\begin{align}
    & E[(M_0(X^n, \bp) - \GT( X^n,\alpha ))^2] \nonumber \\  
    & \ = \ E\left[ \left( \sum_{x \in \cX} p_x I(N_x(X^n) = 0) \right. \right.  \nonumber \\ & \left. \left. \qquad \qquad  - \frac{1}{(1-\alpha)^2 (n-2)}I(N_x(X_2^{n-1}) = 1, E^{\sim x}_{1,n})  \right) ^2 \right] \nonumber \\
     & = \ \sum_{x \in \cX} \bigg( p_x^2 \ \pr(N_x(X^n) = 0)  \nonumber \\ 
     & \qquad + {1}/{((1-\alpha)^4 (n-2)^2)} \ \pr(N_x(X_2^{n-1}) = 1, E^{\sim x}_{1,n}) \bigg) \nonumber \\ 
     &  + \sum_{x \in \cX} \sum_{y \in \cX, y \neq x} p_x p_y \pr(N_x(X^n) = N_y(X^n) = 0) \nonumber \\ 
     & \  - \frac{(1-\alpha)^{-2}}{(n-2)} \ \bigg[p_y \ \pr(N_y(X^n) = 0, N_x(X_2^{n-1}) = 1, E^{\sim x}_{1,n}) \nonumber \\ 
     & \qquad \qquad +  p_x \ \pr(N_x(X^n) = 0, N_y(X_2^{n-1}) = 1, E^{\sim y}_{1,n} ) \bigg] \nonumber \\ 
     & \ + \frac{(1-\alpha)^{-4}}{ (n-2)^2}  \pr(N_x(X_2^{n-1}) = N_y(X_2^{n-1}) = 1, E^{\sim x}_{1,n}, E^{\sim y}_{1,n}) \nonumber \\
     & \  \overset{(a)}{=} \  \sum_{x \in \cX} \bigg( p_x^2 \ Q_x^n(0)  + {1}/{((1-\alpha)^4 (n-2)^2)} \  Q_x^n(1)  \bigg)  \nonumber \\ 
     & \qquad \quad  + \sum_{x \in \cX} \sum_{y \in \cX, y \neq x} T_{x,y}(n,\alpha), \nonumber 
\end{align}
where we get $(a)$ by using the definitions of $Q^n_{x}(a)$, $Q^n_{x,y}(a,b),$ and $T_{x,y}(n,\alpha).$  This completes the proof of \eqref{eq: MSE_expr}.

\subsection{Proof of Lemma \ref{lem: MSE_term_bds}} 
\label{sec: MSE_term_bds_Prf}
Let $x^n$ $\triangleq$ $(x_1,\ldots,x_n),$ $\cX_{\sim x} \triangleq \cX\setminus \{x \}.$ Let $x_{\sim m}$ $\triangleq$ $(x_1,\ldots,x_{m-1},x_{m+1},\ldots,x_n),$ for $m=1$ to $n.$ 
We begin with the proof of \eqref{eq: t2_bd} of Lemma \ref{lem: MSE_term_bds}. 
\subsubsection{Lemma \ref{lem: MSE_term_bds}, Eqn \eqref{eq: t2_bd}}
For any $x \in \cX,$
\begin{align}
     Q_x(0) \ & = \ \pr(N_x(X^n) = 0) \nonumber \\ 
    & = \ \sum_{x^n \in (\cX_{\sim x})^n } \pr( X_i = x_i, 1 \leq i \leq n) \nonumber \\ 
    & \overset{(a)}{=} \ \sum_{x^n \in (\cX_{\sim x})^n } \pr(x_1) \prod_{i=2}^n  \pr( x_i |  x_{i-1})  \nonumber \\ 
    & {=} \ \sum_{x_1 \in \cX\setminus \{x\}} \pr(x_1) \quad \prod_{i=2}^n \left( \sum_{x_i \in \cX\setminus \{x\}} \pr(x_i |  x_{i-1}) \right) \nonumber \\
    & \overset{(b)}{=} \ (1-p_x) \ (1-(1-\alpha)p_x)^{n-1}, \label{eq: Qx_0_expr}  
\end{align}
where we get $(a)$ by using the Markovian property of $X^n,$ $(b)$ by using \eqref{eq: Xn_MC}.
Using \eqref{eq: Qx_0_expr}, we get 
\begin{align}
  \sum_{x \in \cX} p_x^2 \ Q_x(0)  \ & \leq \  \sum_{x \in \cX} p_x^2 \ (1-p_x) \ (1-(1-\alpha)p_x)^{n-1}  \nonumber \\
  & \leq \ \sum_{x \in \cX} p_x^2 \  (1-(1-\alpha)p_x)^{n} \nonumber \\ 
  & \overset{(c)}{\leq} \ \frac{1}{(1-\alpha) (n+1)} \ \sum_{x \in \cX} p_x \ \overset{(d)}{\leq} \ \frac{1}{(1-\alpha) (n-1)}, \nonumber  
\end{align}
where we get $(c)$ by using $\max_{q \in [0,1]} q (1-q)^n \leq 1/(n+1),$ $(d)$ by using $\sum_{x \in \cX} p_x =1$ and $1/(n+1) < 1/(n-1).$ This completes the proof of \eqref{eq: t2_bd}.

\subsubsection{Lemma \ref{lem: MSE_term_bds}, Eqn \eqref{eq: t3_bd}} 
For any $x \in \cX,$
\begin{align}
  Q_x(1) \ & = \ \sum_{m=2}^{n-1} \pr(X_m = x, N_x(X_{\sim m}) = 0) \nonumber   
\end{align}
For $m = 2$ to $n-1,$ 
\begin{align}
  &  \pr(X_m = x, N_x(X_{\sim m}) = 0) \nonumber \\ 
  & = \ \sum_{x_{\sim m} \in (\cX_{\sim x})^{n-1} } \pr(X_m = x, X_l = x_l,l \neq m, 1 \leq l \leq n) \nonumber \\  
  & \overset{(a)}{=} \ \sum_{x_{\sim m} \in (\cX_{\sim x})^{n-1} } \pr(x_1) \ \left(\prod_{i=2}^{m-1}  \pr( x_i |  x_{i-1}) \right) \nonumber \\
  & \qquad \qquad \pr( x| x_{m-1}) \ \pr( x_{m+1}| x)  \  \prod_{i=m+2}^{n}  \pr( x_i |  x_{i-1}) \nonumber \\
  & = \ \sum_{x_1 \in \cX_{\sim x} } \pr( x_1) \ \prod_{i=2}^{m-1} \left( \sum_{x_i \in \cX_{\sim x} } \pr( x_i |  x_{i-1}) \right) \pr( x| x_{m-1}) \nonumber \\
  & \qquad \left( \sum_{x_{m+1} \in \cX_{\sim x} } \pr( x_{m+1}| x) \right)   \prod_{i=m+2}^{n} \left( \sum_{x_i \in \cX_{\sim x} } \pr( x_i |  x_{i-1}) \right) \nonumber \\
  & \overset{(b)}{=} \ (1-p_x) \  (1-(1-\alpha)p_x)^{m-2} \ (1-\alpha) p_x \ (1-\alpha)(1-p_x) \nonumber \\
  & \qquad (1-(1-\alpha)p_x)^{n-m-1} \nonumber \\
  &{=} \ (1-\alpha)^2 \ p_x \ (1-p_x)^2 \ (1-(1-\alpha)p_x)^{n-3}, \label{eq: Qx_1_expr}
\end{align}
where we get $(a)$ by using the Markovian property of $X^n,$ $(b)$ by using \eqref{eq: Xn_MC}.
Therefore
\begin{align}
    &  \sum_{x \in \cX} Q_x(1) \nonumber \\ 
    & = \ {(1-\alpha)^2 (n-2)} \sum_{x \in \cX} p_x \ (1-p_x)^2 \ (1-(1-\alpha)p_x)^{n-3} \nonumber \\
    & \overset{(c)}{\leq} \ {(1-\alpha)^2 (n-2)} \sum_{x \in \cX} p_x \ \leq \ {(1-\alpha)^2 (n-2)}, \nonumber  \end{align}
    where we get $(c)$ by using $(1-p_x)^2 \ (1-(1-\alpha)p_x)^{n-3} \leq 1.$ Dividing by $(1-\alpha)^4 (n-2)^2$ completes the proof of \eqref{eq: t3_bd}.
    
    \subsubsection{Lemma \ref{lem: MSE_term_bds}, Eqn \eqref{eq: t1_bd}} 
    Recall that for $x,y \in \cX,$
    \begin{align}
       & T_{x,y}(n,\alpha) = p_x p_y  Q_{x,y}(0,0) + \frac{1}{ (1-\alpha)^4 (n-2)^2 } Q_{x,y}(1,1) \nonumber \\ 
      & \ - \frac{1}{ (1-\alpha)^2 (n-2) } ( p_x Q_{x,y}(0,1) + p_y Q_{x,y}(1,0) ). \label{eq: Txy_expr}   \end{align}
    Following a method similar to the derivation of \eqref{eq: Qx_0_expr}, we get for any $x,y \in \cX,$ 
    \begin{equation}
        Q_{x,y}(0,0) \  = \ (1-p_x-p_y) \ (1-(1-\alpha)(p_x+p_y))^{n-1}. \label{eq: Qxy_00_expr}
    \end{equation}
    For any $x,y \in \cX,$
\begin{align}
  Q_{x,y}(1,0) \ & = \ \sum_{m=2}^{n-1} \pr(X_m = x, N_x(X_{\sim m}) = N_y(X_{\sim m}) = 0 ). \nonumber   
\end{align}
    Following a method similar to the derivation of \eqref{eq: Qx_1_expr}, we get for any $x,y \in \cX,$ 
    \begin{align}
        Q_{x,y}(1,0) \ & = \ (n-2) \ (1-\alpha)^2 \ p_x \ (1-p_x-p_y)^2 \nonumber \\ 
        &  \qquad (1-(1-\alpha)(p_x+p_y))^{n-3} \label{eq: Qxy_10_expr}
    \end{align}
For $m_1 = 2$ to $n-2,$ $m_2 = m_1+1$ to $n-1,$  let $X_{\sim m_1, m_2} \triangleq (X_1,\ldots, X_{m_1-1}, X_{m_1+1}, \ldots, X_{m_2-1},X_{m_2+1},\ldots  X_n)$ denote the samples from $X_1$ to $X_n$ except $X_{m_1}, X_{m_2}.$ To reduce clutter, we define the event $E^{x,y}_{\sim m_1, m_2} \triangleq  ( N_x(X_{\sim m_1, m_2}) = N_y(X_{\sim m_1,m_2}) = 0 )$.
  For any $x,y \in \cX,$
\begin{align}
   Q_{x,y}(1,1) 
  & =  \sum_{m_1=2}^{n-2} \sum_{m_2=m_1+1}^{n-1} \pr \left(X_{m_1} = x, X_{m_2} = y, E^{x,y}_{\sim m_1, m_2}\right) \nonumber  \\
  & \qquad \qquad \quad  + \pr \left(X_{m_1} = y, X_{m_2} = x, E^{x,y}_{\sim m_1, m_2} \right). \label{eq: Qxy_11_interim}
\end{align}
 Following a method similar to the derivation of \eqref{eq: Qx_1_expr}, we get for any $x,y \in \cX,$ 
 \begin{align}
  & \pr\left(X_{m_1} = x, X_{m_1+1} = y, E^{x,y}_{\sim m_1, m_1+1}\right)  \nonumber \\ 
  & \ =  (1-\alpha)^3 \ p_x p_y \ (1-p_x-p_y)^2 \ (1-(1-\alpha)(p_x+p_y))^{n-4} \nonumber \\ 
  & \qquad \qquad \qquad \text{for } m_1 = 2 \text{ to }  n-2, m_2 = m_1+1, \label{eq: Qxy_11_consec} \\ 
  & \pr\left(X_{m_1} = x, X_{m_2} = y, E^{x,y}_{\sim m_1, m_2}\right)  \nonumber \\ 
  & \ =  (1-\alpha)^4 \ p_x p_y \ (1-p_x-p_y)^3 \ (1-(1-\alpha)(p_x+p_y))^{n-5} \nonumber \\ 
  & \qquad \qquad \qquad \text{for } m_1 = 2 \text{ to } n-3, m_2 \geq m_1+2. \label{eq: Qxy_11_nonconsec} 
 \end{align}
 Using \eqref{eq: Qxy_11_consec}, \eqref{eq: Qxy_11_nonconsec} in \eqref{eq: Qxy_11_interim}, we get 
 \begin{align}
      Q_{x,y}(1,1)  
   & = \ (n-3) \ (1-\alpha)^3 \ p_x p_y \ (1-(1-\alpha)(p_x+p_y))^{n-5}  \nonumber \\ & \qquad (1-p_x-p_y)^2  \ \left[2(1-(1-\alpha)(p_x+p_y)) \right. \nonumber \\ 
   & \left. \qquad \qquad \qquad +  (n-4)(1-\alpha) \ (1-p_x-p_y) \right], \label{eq: Qxy_11_expr}
 \end{align}
 Substituting \eqref{eq: Qxy_00_expr}, \eqref{eq: Qxy_10_expr}, and \eqref{eq: Qxy_11_expr} in \eqref{eq: Txy_expr}, simplifying and taking the summation over $\{(x,y): x,y \in \cX, x \neq y\}$, we get
 \begin{align}
   & \sum_{x,y \in \cX: x \neq y} T_{x,y}(n,\alpha) \nonumber \\ 
    & = \ \sum_{x,y \in \cX: x \neq y} p_x p_y \ (1-p_x-p_y) \ (1-(1-\alpha)(p_x+p_y))^{n-5} \nonumber \\ 
    & \qquad  \left( (p_x+p_y)^2 \ \left[ \alpha^2 - (1-\alpha)^2(1-p_x-p_y) \right]^2 \right. \nonumber \\ 
    & \qquad  \left. +  (1-p_x-p_y)  \left[ 2 \frac{(n-3)(1-\alpha)^{-1}}{(n-2)^2} (1-(1-\alpha)(p_x+p_y)) \right. \right. \nonumber \\ & \left. \left. \qquad \qquad \qquad \qquad \ - \left(\frac{3}{(n-2)} - \frac{2}{(n-2)^2}\right) (1-p_x-p_y)  \right] \right) \nonumber \\
    & \overset{(a)}{\leq} \ \sum_{x,y \in \cX: x \neq y} p_x p_y \left( \left( \alpha^2 \left( \frac{\alpha^2}{(1-\alpha)^2} + 2 \right) +  (1-\alpha)^2 \right) O(1/n^2) \right. \nonumber  \\ 
    & \ \qquad \qquad \qquad \qquad \qquad + \left. \left[ 2 \frac{(1-\alpha)^{-1}}{(n-2)} + O(1/n) \right] \right), \nonumber \\ 
    & \overset{(b)}{\leq} \ \left( \alpha^2 \left( \frac{\alpha^2}{(1-\alpha)^2} + 2 \right) +  (1-\alpha)^2 \right) O(1/n^2) \nonumber \\ 
    & \ \qquad \qquad \qquad  +  2 \frac{(1-\alpha)^{-1}}{(n-2)} + O(1/n), \nonumber
 \end{align}
where we get $(a)$ using $\max_{p \in [0,1]} p^2(1-p)^m \leq 4/m^2$ and $ (1-p_x-p_y) \ (1-(1-\alpha)(p_x+p_y)) \leq 1 ,$  $(b)$ by using $ \sum_{x,y \in \cX: x \neq y} p_x p_y \leq \sum_{x,y \in \cX} \ p_x p_y = 1. $
This completes the proof of \eqref{eq: t1_bd}.

\subsection{Proof of Lemma \ref{lem:Lecam_lb} }
\label{subsec:Lecam_lb_Prf}

\begin{align}
  \text{ Let }    & M \  \triangleq \ 2 \  \sup_{Q \in \cQ} E [ d( \widehat{\theta}(Y) , \theta(Y,Q) ) ] \nonumber \\ 
    & \geq \ E_{Q_1} [ d( \widehat{\theta}(Y) , \theta(Y,Q_1) ) ] \ + \ E_{Q_2} [ d( \widehat{\theta}(Y) , \theta(Y,Q_2) ) ] \nonumber \\ 
    &  \overset{(a)}{\geq} \ E_{Q_1} [ d(\widehat{\theta}(Y), \cD_1) ] \ + \ E_{Q_2} [ d(\widehat{\theta}(Y), \cD_2) ] - (\epsilon_1 + \epsilon_2) \Delta \label{eq:Lecam_bd_int}
\end{align}
where $E_{Q_i}[\ \cdot \ ]$ denotes expectation with respect to the distribution $Q_i,$  $d(\widehat{\theta}(Y), \cD_i) \triangleq \min_{u \in \cD_i} d(\widehat{\theta}(Y) , u) ,$ $i = 1,2,$ and we get $(a)$ by using the following lemma.

\begin{lemma}
Let $ g_1$ and $g_2$ be two non-negative functions of a random variable $Z$ such that $g_1(Z) \geq g_2(Z)$ with probability atleast $1 - \epsilon.$ Then 
\begin{align}
    E[g_1(Z)] \ & \geq E[g_2(Z)] - \epsilon \ \sup_{z} g_2(z), \label{eq:Exp_lb} 
\end{align}
\label{lem:Exp_lb} 
\end{lemma}
\begin{proof}
To get \eqref{eq:Exp_lb}, we first begin with $E[g_2(Z)].$ 
\begin{align}
     E[g_2(Z)] \
    &  = \ E[g_2(Z) \ I(g_2(Z) \leq g_1(Z))] \nonumber \\ 
    & \qquad  + \ E[g_2(Z) \ I(g_2(Z) > g_1(Z))] \nonumber \\ 
    & \leq  E[g_2(Z) \ I(g_2(Z) \leq g_1(Z))] \nonumber \\ 
    & \qquad \qquad + \left( \sup_{z} \  g_2(z) \right)  \pr( g_2(Z) > g_1(Z)  ) \nonumber \\ 
    & \overset{(a)}{\leq} E[g_2(Z) \ I(g_2(Z) \leq g_1(Z))] \ + \ \epsilon \  \sup_{z} \  g_2(z)   \label{eq:Exp_lb_int}
\end{align}
where we get $(a)$ by using $\pr(g_1(Z) \geq g_2(Z)) \geq 1-\epsilon.$ 
\begin{align}
    E[g_1(Z)] \ & = E[g_1(Z) \ I(g_2(Z) \leq g_1(Z))] \nonumber \\ 
    & \qquad + \ E[g_1(Z) \ I(g_2(Z) > g_1(Z))] \nonumber  \\
    & \geq E[g_1(Z) \ I(g_2(Z) \leq g_1(Z))] \nonumber \\ 
    & \geq E[g_2(Z) \ I(g_2(Z) \leq g_1(Z))] \nonumber \\ 
    &  \overset{(b)}{\geq}   E[g_2(Z)] - \epsilon \   \sup_{z} \  g_2(z) \nonumber 
\end{align}
where we get $(b)$ by using \eqref{eq:Exp_lb_int}.
\end{proof}
Setting $g_1(Y)$ to $d( \widehat{\theta}(Y) , \theta(Y,Q_i) ),$ $g_2(Y)$ to $ d(\widehat{\theta}(Y), \cD_i),$ $i = 1,2,$ and applying the above lemma gives \eqref{eq:Lecam_bd_int}. Note that $\Delta = \max_{u,v \in \cD} d(u,v) \geq  \ d(\widehat{\theta}(Y), \cD_i),$ $i = 1,2.$ \\ 

Since $ d(\widehat{\theta}(Y), \cD_1) + d(\widehat{\theta}(Y), \cD_2) \geq d(\cD_1, \cD_2) = 2\delta,$ we get $ \frac{1}{2\delta} \ ( d(\widehat{\theta}(Y), \cD_1) +  d(\widehat{\theta}(Y), \cD_2) ) \geq 1.$ Therefore 
\begin{align}
    M \ & \geq \ 2\delta \left( E_{Q_1} \left[ \frac{1}{2\delta} \ d(\widehat{\theta}(Y), \cD_1) \right]  +  E_{Q_2} \left[ \frac{1}{2\delta} \ d(\widehat{\theta}(Y), \cD_2) \right] \right) \nonumber \\ 
    & \qquad \qquad - (\epsilon_1 + \epsilon_2) \Delta  \nonumber \\
    & \geq \ 2\delta \left( \inf_{h_1, h_2 \geq 0, h_1 + h_2 =1} E_{Q_1} [h_1(Y)] + E_{Q_2}  [h_2(Y)] \right) \nonumber \\ 
    & \qquad \qquad - (\epsilon_1 + \epsilon_2) \Delta, \nonumber 
\end{align}
where $h_1$ and $h_2$ are arbitrary functions of $Y$. Using the standard fact   
$$\inf_{h_1, h_2 \geq 0, h_1 + h_2 =1} E_{Q_1} [h_1(Y)] + E_{Q_2}  [h_2(Y)] = \lVert Q_1 \wedge Q_2\rVert,$$ 
completes the proof. 

\subsection{Proof of Lemma \ref{lem:TV_bd} }
\label{subsec:TV_bd_Prf}
To reduce clutter, we denote $P_1\triangleq P_{\bp(0,L)}$, $P_2\triangleq P_{\bp(\beta,L)}$, $\bp_1=\bp(0,L)$ and $\bp_2=\bp(\beta,L)$.

To show the bound in \eqref{eq:TV_bd}, we first bound the total variation distance between $P_1(X^n)$ and $P_2(X^n)$ by the KL divergence $D_{KL}(P_1(X^n) || P_2(X^n))$ using Pinsker's inequality. 
\begin{lemma}
\textbf{Pinkser's inequality} \cite[Theorem 4.19]{boucheron2013}
\begin{equation}
      ||P_1(X^n) - P_2(X^n)||_{TV} \ \leq \ \frac{1}{\sqrt{2}} \ \sqrt{ D_{KL}(P_1(X^n) || P_2(X^n)) } \label{eq:Pinsker's_ineq}
\end{equation}
where $ D_{KL}(P_1(X^n) || P_2(X^n)) $ is the KL divergence between $P_1(X^n)$ and $P_2(X^n).$
\end{lemma}
Our next lemma expresses $ D_{KL}(P_1(X^n) || P_2(X^n)) $ in terms of the KL divergence between the  distributions $\bp_1, \bp_2$ and the KL divergence between the corresponding rows of $P_1$ and $P_2.$
\begin{lemma}
For any two t.p.ms $P_1, P_2,$ on an $\cX,$ with $\bp_1, \bp_2$ as their respective stationary distributions, 
\begin{align}
    & D_{KL}(P_1(X^n) || P_2(X^n)) \nonumber \\
    &  =  D_{KL}(\bp_1 || \bp_2) + (n-1) \ \sum_{x \in \cX} \pi_{1,x} \ D_{KL}( P_1(\cdot | x) \ || \ P_2(\cdot | x) ) \label{eq:KL_div_Markov}
\end{align}
where $P_i( \cdot |x)$ denotes the row of the t.p.m $P_i, i = 1,2,$  with transition probabilities from the state $x.$
\label{lem:KL_div_Markov}
\end{lemma}
\begin{proof}
Let $x^n \triangleq (x_1, x_2,\ldots,x_n) \in \cX^n.$
\begin{align}
    & D_{KL}(P_1(X^n) || P_2(X^n)) \nonumber \\
    & = \ \sum_{x^n \in \cX^n} P_1(x^n) \ \ln (P_1(x^n)/P_2(x^n)) \nonumber \\ 
    & \overset{(a)}{=} \ \sum_{x^n \in \cX^n}  \ P_1( x^n)  \bigg[ \ln \ (\pi_{1,x_1}/\pi_{2,x_1})  \nonumber \\ 
    &  \quad \qquad \qquad \qquad \qquad + \sum_{l = 2}^n \ \ln  \frac{P_1(X_l = x_l | X_{l-1} = x_{l-1})}{ P_2(X_l = x_l | X_{l-1} = x_{l-1}) }  \bigg] \nonumber \\
    & \overset{(b)}{=} D_{KL}(\bp_1||\bp_2) + \sum_{l=2}^n \sum_{x^l \in \cX^l} P_1(x^l) \ln \frac{P_1(X_l = x_l | X_{l-1} = x_{l-1})}{ P_2(X_l = x_l | X_{l-1} = x_{l-1}) } \nonumber \\
    & \overset{(c)}{=} D_{KL}(\bp_1||\bp_2) + \sum_{l=2}^n \sum_{x_{l-1}, x_l \in \cX}  p_{1,x_{l-1}}   P_1(X_l = x_l | X_{l-1} = x_{l-1}) \nonumber \\ 
    & \qquad \qquad \qquad \qquad \qquad   \ln \frac{P_1(X_l = x_l | X_{l-1} = x_{l-1})}{ P_2(X_l = x_l | X_{l-1} = x_{l-1}) } \nonumber \\ 
    & = D_{KL}(\bp_1||\bp_2) + (n-1) \sum_{ x \in \cX} p_{1,x} D_{KL}( P_1(\cdot | x) || P_2(\cdot | x) ), \nonumber 
\end{align}
where we get $(a)$ by using the Markov property $ P_i(x^n) = p_{i,x_1} \ \prod_{l-2}^n P_i(X_l = x_l| X_{l-1} = x_{l-1}),$ $i = 1,2,$ $(b)$ and $(c)$ by appropriately marginalizing $P_1(x^n).$
\end{proof}
Using the values specified for $p_{i,x},$  $P_i(X_2 = y|X_1 = x)$ for $x,y \in \{1,\ldots, L+1 \},$ $i = 1,2,$ in Section \ref{subsec:Lecam_lb}, we get
\begin{align}
    &D_{KL}(\bp_1||\bp_2) = -0.5 \ln \left( { 1 - {4 \gamma^2} } \right) \label{eq:p_DKL} \\
    &D_{KL}(P_1(\cdot|1) || P_2(\cdot|1) ) = -0.5 \ (1 + \alpha)  \ln  \left( {{1+2\gamma\frac{(1-\alpha)}{(1+\alpha)}}}\right) \nonumber \\ 
    & \qquad\qquad\qquad\qquad - 0.5 \ ( 1 - \alpha) \ \ln \left( 1 - 2\gamma \right) \nonumber \\
    &\qquad\leq \ - 0.5 (1-\alpha) \ln (1-4\gamma^2). \label{eq:p_DKL_r1}
\end{align}
For $x \in \{ 2,\ldots, L+1 \}$,
\begin{align}
    &D_{KL}(P_1(\cdot|x) || P_2(\cdot|x) )=-(1-\alpha)\nonumber\\
    &\qquad \bigg( 0.5 \ln (1{-}4\gamma^2) + {\gamma^2}\left[ \frac{ 4 \alpha  (1{-}2 \gamma)^{-1}  }{ (2L \alpha + (1{-}\alpha) (1{-}2 \gamma) ) }  \right] \bigg),  \nonumber \\
    & \leq - 0.5 (1-\alpha) \ln (1-4\gamma^2). \label{eq:p_DKL_r2}
\end{align}
Using \eqref{eq:p_DKL}, \eqref{eq:p_DKL_r1}, and \eqref{eq:p_DKL_r2} in \eqref{eq:KL_div_Markov}, we get 
\begin{align}
     & D_{KL}(P_1(X^n) || P_2(X^n)) \nonumber\\
     &\quad \leq  -0.5  [1 + (n-1)(1-\alpha) ]  \ln (1-4\gamma^2)  \nonumber \\
     &\quad \overset{(a)}{\leq} 4 [1 + (n-1)(1-\alpha) ] \gamma^2, \nonumber 
\end{align}
where we get $(a)$ by using $-\ln(1-x) \leq 2x,$ for $x \in (0,0.5)$. Plugging the above bound into \eqref{eq:Pinsker's_ineq} completes the proof of Lemma \ref{lem:TV_bd}.

\bibliographystyle{IEEEtran}
\bibliography{refs}

\begin{thebibliography}{10}
\providecommand{\url}[1]{#1}
\csname url@samestyle\endcsname
\providecommand{\newblock}{\relax}
\providecommand{\bibinfo}[2]{#2}
\providecommand{\BIBentrySTDinterwordspacing}{\spaceskip=0pt\relax}
\providecommand{\BIBentryALTinterwordstretchfactor}{4}
\providecommand{\BIBentryALTinterwordspacing}{\spaceskip=\fontdimen2\font plus
\BIBentryALTinterwordstretchfactor\fontdimen3\font minus
  \fontdimen4\font\relax}
\providecommand{\BIBforeignlanguage}[2]{{%
\expandafter\ifx\csname l@#1\endcsname\relax
\typeout{** WARNING: IEEEtran.bst: No hyphenation pattern has been}%
\typeout{** loaded for the language `#1'. Using the pattern for}%
\typeout{** the default language instead.}%
\else
\language=\csname l@#1\endcsname
\fi
#2}}
\providecommand{\BIBdecl}{\relax}
\BIBdecl

\bibitem{Chen1996}
S.~F. Chen and J.~Goodman, ``An empirical study of smoothing techniques for
  language modeling,'' in \emph{Proceedings of the 34th Annual Meeting on
  Association for Computational Linguistics}, ser. ACL '96, 1996, pp. 310--318.

\bibitem{Good1953}
I.~J. Good, ``The population frequencies of species and the estimation of
  population parameters,'' \emph{Biometrika}, vol.~40, no. 3/4, pp. 237--264,
  1953.

\bibitem{Gale1995}
W.~A. Gale and G.~Sampson, ``Good‐{T}uring frequency estimation without
  tears,'' \emph{Journal of Quantitative Linguistics}, vol.~2, no.~3, pp.
  217--237, 1995.

\bibitem{Gale1991}
K.~W.Church and W.~A.Gale, ``Probability scoring for spelling correction,''
  \emph{Statistics and Computing}, vol.~1, pp. 93--103, 1991.

\bibitem{Song99}
F.~Song and W.~B. Croft, ``A general language model for information
  retrieval,'' ser. CIKM '99.\hskip 1em plus 0.5em minus 0.4em\relax
  Association for Computing Machinery, 1999, p. 316–321.

\bibitem{Chao1992}
A.~Chao and S.-M. Lee, ``Estimating the number of classes via sample
  coverage,'' \emph{Journal of the American Statistical Association}, vol.~87,
  no. 417, pp. 210--217, 1992.

\bibitem{Colwell12}
R.~K. Colwell, A.~Chao, N.~J. Gotelli, S.-Y. Lin, C.~X. Mao, R.~L. Chazdon, and
  J.~T. Longino, ``Models and estimators linking individual-based and
  sample-based rarefaction, extrapolation and comparison of assemblages,''
  \emph{Journal of Plant Ecology}, vol.~5, no.~1, pp. 3--21, 2012.

\bibitem{McAllester2000}
D.~A. McAllester and R.~E. Schapire, ``On the convergence rate of
  {G}ood-{T}uring estimators,'' in \emph{Proceedings of the Thirteenth Annual
  Conference on Computational Learning Theory}, 2000, pp. 1--6.

\bibitem{Chandra19}
P.~Chandra, A.~Pradeep, and A.~Thangaraj, ``Improved tail bounds for missing
  mass and confidence intervals for {G}ood-{T}uring estimator,'' in
  \emph{National Conference on Communications 2019}, Feb 2019, p. to appear.

\bibitem{McAllester2003}
D.~McAllester and L.~Ortiz, ``Concentration inequalities for the missing mass
  and for histogram rule error,'' \emph{J. Mach. Learn. Res.}, vol.~4, pp.
  895--911, Dec. 2003.

\bibitem{berend2013}
D.~Berend and A.~Kontorovich, ``On the concentration of the missing mass,''
  \emph{Electron. Commun. Probab.}, vol.~18, p. 7 pp., 2013.

\bibitem{ben-hamou2017}
A.~Ben-Hamou, S.~Boucheron, and M.~I. Ohannessian, ``Concentration inequalities
  in the infinite urn scheme for occupancy counts and the missing mass, with
  applications,'' \emph{Bernoulli}, vol.~23, no.~1, pp. 249--287, 02 2017.

\bibitem{Chandra19_ISIT}
P.~{Chandra} and A.~{Thangaraj}, ``Concentration and tail bounds for missing
  mass,'' in \emph{2019 IEEE International Symposium on Information Theory
  (ISIT)}, July 2019, pp. 1862--1866.

\bibitem{Rajaraman17}
N.~Rajaraman, A.~Thangaraj, and A.~T. Suresh, ``Minimax risk for missing mass
  estimation,'' in \emph{2017 IEEE International Symposium on Information
  Theory (ISIT)}, June 2017, pp. 3025--3029.

\bibitem{Acharya17}
J.~{Acharya}, Y.~{Bao}, Y.~{Kang}, and Z.~{Sun}, ``Improved bounds for minimax
  risk of estimating missing mass,'' in \emph{2018 IEEE International Symposium
  on Information Theory (ISIT)}, June 2018, pp. 326--330.

\bibitem{Acharya13}
J.~Acharya, A.~Jafarpour, A.~Orlitsky, and A.~T. Suresh, ``Optimal probability
  estimation with applications to prediction and classification,'' in
  \emph{Proceedings of the 26th Annual Conference on Learning Theory}, vol.~30,
  2013, pp. 764--796.

\bibitem{Orlitsky2015}
A.~Orlitsky and A.~T. Suresh, ``Competitive distribution estimation: Why is
  {G}ood-{T}uring good,'' in \emph{Advances in Neural Information Processing
  Systems 28}, 2015, pp. 2143--2151.

\bibitem{Abs_disc_17}
M.~Falahatgar, M.~I. Ohannessian, A.~Orlitsky, and V.~Pichapati, ``The power of
  absolute discounting: all-dimensional distribution estimation,'' in
  \emph{Advances in Neural Information Processing Systems}, vol.~30, 2017.

\bibitem{Ohannessian19}
E.~Mossel and M.~I. Ohannessian, ``On the impossibility of learning the missing
  mass,'' \emph{Entropy}, vol.~21, no.~1, 2019.

\bibitem{hao_doubly_comp19}
Y.~Hao and A.~Orlitsky, ``Doubly-competitive distribution estimation,'' in
  \emph{Proceedings of the 36th International Conference on Machine Learning},
  ser. Proceedings of Machine Learning Research, vol.~97.\hskip 1em plus 0.5em
  minus 0.4em\relax PMLR, 2019, pp. 2614--2623.

\bibitem{Skorski21}
M.~Skorski, ``Mean-squared accuracy of good-turing estimator,'' in \emph{2021
  IEEE International Symposium on Information Theory (ISIT)}, 2021, pp.
  2846--2851.

\bibitem{Painsky21}
A.~Painsky, ``Refined convergence rates of the good-turing estimator,'' in
  \emph{2021 IEEE Information Theory Workshop (ITW)}, 2021, pp. 1--5.

\bibitem{Olgica_15}
S.~M. H.~T. Yazdi, H.~M. Kiah, E.~Garcia-Ruiz, J.~Ma, H.~Zhao, and
  O.~Milenkovic, ``Dna-based storage: Trends and methods,'' \emph{IEEE
  Transactions on Molecular, Biological and Multi-Scale Communications},
  vol.~1, no.~3, pp. 230--248, 2015.

\bibitem{Olgica_21}
E.~Chien, O.~Milenkovic, and A.~Nedich, ``Support estimation with sampling
  artifacts and errors,'' in \emph{2021 IEEE International Symposium on
  Information Theory (ISIT)}, 2021, pp. 244--249.

\bibitem{Olgica_09}
F.~Farnoud, O.~Milenkovic, and N.~P. Santhanam, ``Small-sample distribution
  estimation over sticky channels,'' in \emph{2009 IEEE International Symposium
  on Information Theory}.\hskip 1em plus 0.5em minus 0.4em\relax IEEE, 2009,
  pp. 1125--1129.

\bibitem{Olgica_12}
F.~Farnoud, N.~P. Santhanam, and O.~Milenkovic, ``Alternating markov chains for
  distribution estimation in the presence of errors,'' in \emph{2012 IEEE
  International Symposium on Information Theory Proceedings}.\hskip 1em plus
  0.5em minus 0.4em\relax IEEE, 2012, pp. 2017--2021.

\bibitem{nanopore_16}
A.~Magner, J.~Duda, W.~Szpankowski, and A.~Grama, ``Fundamental bounds for
  sequence reconstruction from nanopore sequencers,'' \emph{IEEE Transactions
  on Molecular, Biological and Multi-Scale Communications}, vol.~2, no.~1, pp.
  92--106, 2016.

\bibitem{fried2021alphalazy}
S.~Fried and G.~Wolfer, ``On the $\alpha$-lazy version of markov chains in
  estimation and testing problems,'' 2021.

\bibitem{markov_book}
R.~G. Gallager, ``Finite state {M}arkov chains,'' in \emph{Discrete Stochastic
  Processes}.\hskip 1em plus 0.5em minus 0.4em\relax Springer, 1996, pp.
  103--147.

\bibitem{yu1997lecam}
B.~Yu, ``Assouad, fano, and le cam,'' in \emph{Festschrift for Lucien Le
  Cam}.\hskip 1em plus 0.5em minus 0.4em\relax Springer, 1997, pp. 423--435.

\bibitem{boucheron2013}
S.~Boucheron, G.~Lugosi, and P.~Massart, \emph{Concentration Inequalities: A
  Nonasymptotic Theory of Independence}.\hskip 1em plus 0.5em minus 0.4em\relax
  Oxford, UK: Oxford University Press, 2013.

\end{thebibliography}

\end{document}